\documentclass[11pt]{article}
\usepackage{amsmath,latexsym}
\usepackage{amssymb}

\usepackage{theorem}
\numberwithin{equation}{section}

\title{Curvature and Uniformization}
\author{Rafe Mazzeo \thanks{Research supported in part by NSF Grant 
DMS-9971975 and also at MSRI by NSF grant DMS-9701755. Email: 
mazzeo\@@math.stanford.edu} \\ Stanford University 
\and Michael Taylor \thanks{Research supported in part by
NSF Grant DMS-9877077. Email: met\@@math.unc.edu} \\ University of North Carolina}
\date{}

\newcommand{\bW}{\partial\Omega}
\newcommand{\Om}{\Omega}
\newcommand{\Ombar}{\overline{\Omega}}

\newcommand{\pa}{\partial}
\newcommand{\del}{\partial}
\newcommand{\RR}{\mathbb R}
\newcommand{\CC}{\mathbb C}
\newcommand{\TT}{\mathbb T}
\newcommand{\ZZ}{\mathbb Z}
\newcommand{\calA}{\mathcal A}
\newcommand{\calC}{\mathcal C}

\newcommand{\tv}{\tilde{v}}
\newcommand{\tQ}{\tilde{Q}}

\newcommand{\ep}{\varepsilon}
\newcommand{\beq}{\begin{equation}}
\newcommand{\eeq}{\end{equation}}
\def\beqa{\begin{eqnarray}}
\def\eeqa{\end{eqnarray}}

\theoremstyle{plain}
\newtheorem{theorem}{Theorem}[section]
\newtheorem{proposition}[theorem]{Proposition}
\newtheorem{corollary}[theorem]{Corollary}
\newtheorem{lemma}[theorem]{Lemma}

\begin{document}

\maketitle

\begin{abstract}
We approach the problem of uniformization of general Riemann
surfaces through consideration of the curvature equation,
and in particular the problem of constructing Poincar\'e metrics
(i.e., complete metrics of constant negative curvature) by
solving the equation $\Delta u - e^{2u} = K_0(z)$ on general
open surfaces. A few other topics are discussed, including
boundary behavior of the conformal factor $e^{2u}$ giving the
Poincar\'e metric when the Riemann surface has smoothly
bounded compact closure, and also a curvature equation
proof of Koebe's disk theorem.
\end{abstract}

\section{Introduction}

Let $M$ be a smooth, connected, oriented two-dimensional manifold. 
A Riemannian metric $g$ on $M$ determines a conformal class
\[
[g] = \{e^{2u}g: u \in C^\infty(M)\},
\]
and there is a well-known bijection between the set of conformal classes 
and the set of complex structures on $M$. A Riemann surface is such
a surface endowed with a particular choice of conformal or complex structure.

It is reasonable to seek a canonical metric in each conformal class and a 
natural candidate is one with constant Gaussian curvature $K$. The case 
where $K$ is negative arises most frequently, and in any case is the one
upon which we mostly concentrate. Thus we define a 
{\it Poincar\'e metric} on a Riemann surface $M$ to be one 
(in the conformal class of $M$) that is complete and that has 
Gauss curvature $K \equiv -1$. A basic example is the Poincar\'e metric 
$G$ on the unit disk $D_1 \subset \RR^2$, which has components
\begin{equation}
G_{jk}=\frac{4}{(1-r^2)^2}\, \delta_{jk},\quad r^2=x_1^2+x_2^2.
\label{1.1}
\end{equation}
This is the unique Poincar\'e metric in $[\delta]$, and it is invariant
with respect to all conformal (or holomorphic) automorphisms of $D_1$. 

If $g_0$ is a metric on $M$, with Gauss curvature function $K_0(x)$, then 
$g = e^{2u}g_0$ has Gauss curvature $K=(K_0-\Delta u)e^{-2u}$, which is 
equal to $-1$ provided $u$ satisfies 
\begin{equation}
\Delta u - e^{2u} = K_0(x).
\label{1.2}
\end{equation}
In particular, to find a Poincar\'e metric $g \in [g_0]$ it is sufficient
to solve (\ref{1.2}) and show that $e^{2u}g_0$ is complete.

Poincar\'e metrics are closely related to conformal (holomorphic) coverings 
by $D_1$. In fact, a Poincar\'e metric $g$ on $M$ lifts to a Poincar\'e 
metric $\tilde{g}$ on the universal cover $\widetilde{M}$, and the
covering map $\varphi: \widetilde{M} \to M$ is by definition a local isometry,
hence conformal. On the other hand, a basic theorem in differential geometry 
asserts that $(\widetilde{M},\widetilde{g})$ is isometric to the disk $D_1$
with its Poincar\'e metric (\ref{1.1}). Therefore we obtain a 
holomorphic covering map
\[
\varphi: D_1 \longrightarrow M,
\]
which is a local isometry between the Poincar\'e metrics on $D_1$ and $M$.
Conversely, if $\varphi$ is any such conformal covering map, 
the deck transformations on $D_1$ are conformal and thus fix the
Poincar\'e metric there. Hence $\varphi$ acts by isometries and pushes 
down to a Poincar\'e
metric on $M$. Extending this reasoning slightly, we see that
\begin{equation}
\mbox{{\it If a Poincar\'e metric exists on a Riemann surface $M$, it 
is unique.}}
\label{1.3}
\end{equation}

The discussion in the last paragraph makes it clear that
the construction of Poincar{\'e} metrics is intimately related to
the classical uniformization theorem, which we now state:

\smallskip
\noindent \textbf{Uniformization Theorem} 
{\it Every simply-connected Riemann surface is 
holomorphically equivalent to either $\widehat{\CC}$, $\CC$, or $D_1$.}
\smallskip

\noindent Here $\widehat{\CC}$ denotes the Riemann sphere.
An equivalent statement is that any (connected) Riemann surface $M$ can be
holomorphically covered by $\widehat{\CC},\ \CC$, or $D_1$.
It is well known that this result can be established when $M$ is compact by 
directly solving the curvature equation; cf.~\S{8} for further discussion
of this.  One of our goals here is to give a direct treatment of the
curvature equation on a broad class of Riemann surfaces, and to use this 
to establish the uniformization theorem.

We proceed in a series of relatively easy steps. In \S{2} we 
commence by finding a Poincar\'e metric when $M$ is the interior of a 
compact smooth surface with boundary. Section 3 takes up another theme,
the boundary behavior of the Poincar\'e metric in this case.
In \S 4, an approximation argument is used to 
produce a Poincar\'e metric on any domain $\Omega$ in the complex plane
whose complement has at least two points.
In \S{5} we take the space to advertise a purely curvature proof
of Koebe's disk theorem, and its well known corollary 
about normality of a family of univalent maps.
Section 6 establishes the uniformization theorem for general
simply connected Riemann surfaces, as a consequence of results of 
sections 2 and 5.  In \S{7} we relate the  dichotomy
between Riemann surfaces covered by $D_1$ and those covered
by $\CC$ to a dichotomy in the behavior of the curvature equation. 
In \S{8} we discuss the uniformization theorem for 
compact surfaces.

We say more about the second main theme of this paper, taken up in \S{3}.
Many developments in modern function theory have focused on the 
connection between the regularity of the boundary of $M$ (especially 
when it is a planar domain) and
the regularity of the mapping $\varphi$.  From the point of view
here, it seems also of interest to examine the boundary behavior of
the solution to (\ref{1.2}) yielding
the Poincar{\'e} metric, especially when $M$ has compact
closure in a larger Riemann surface.  The boundary regularity results
we obtain in \S{3} on $e^{-u}$ have implications for the qualitative
behavior of the covering map $D_1\rightarrow M$.

We conclude this introduction by providing a few explicit examples of 
Poincar\'e metrics to illustrate
various phenomena that can occur.  Also, we will have specific
use for several of these formulas in \S{3} and \S{4}.
\begin{itemize}
\item The upper half-plane $H^+=\{x\in\RR^2: x_2>0\}$ has Poincar\'e metric 
\begin{equation}
g_{jk}=x^{-2}_2\, \delta_{jk}.
\label{1.4}
\end{equation}
This may be obtained from (\ref{1.1}) using the standard linear 
fractional transformation that maps $D_1$ to $H^+$.  
\item The Poincar\'e metric on the punctured disk 
$D^*=\{x\in\RR^2:0<|x|<1\}$ is
\begin{equation}
g_{jk}=\Bigl(r\, \log \frac{1}{r}\Bigr)^{-2} \delta_{jk},
\label{1.5}
\end{equation}
as can be verified using the covering $H^+\rightarrow D^*$, 
$z \mapsto e^{iz}$. 
\item
The strip $\Sigma=\{x\in\RR^2:0<x_2<\pi\}$ has Poincar{\'e} metric
\begin{equation}
g_{jk}=(\sin x_2)^{-2}\, \delta_{jk},
\label{1.6}
\end{equation}
as one obtains from (\ref{1.4}) via the conformal diffeomorphism
$\Sigma\rightarrow H^+,\ z\mapsto e^z$.
\item
The annulus $A_b=\{x\in\RR^2:e^{-\pi/b}<|x|<1\}$ has Poincar{\'e} metric
\begin{equation}
g_{jk}=\Bigl[ \frac{b}{r\sin(b\,\log \frac{1}{r})}\Bigr]^2 \delta_{jk},
\label{1.7}
\end{equation}
as can be seen using (\ref{1.6}) and the covering $\Sigma\rightarrow A_b,
\ z\mapsto e^{iz/b}$.  Note that the $b\rightarrow 0$ limit gives
(\ref{1.5}).
\item The quarter-plane $Q = \{x \in\RR^2: x_1 > 0, x_2 > 0\}$ has
Poincar\'e metric 
\begin{equation}
g_{jk} = \frac{x_1^2+x_2^2}{x_1^2 x_2^2}\, \delta_{jk},
\label{1.8}
\end{equation}
as one obtains from (\ref{1.4}) via the map $Q\rightarrow H^+,\ z\mapsto
z^2$.
\end{itemize}

\section{Smoothly bounded Riemann surfaces}

Let $\Ombar$ be a compact, oriented, connected 2-dimensional Riemannian 
manifold with smooth boundary, with metric $g_0$.  We can suppose that 
$\Ombar$ is contained in a larger open Riemann surface $M$. We produce
a solution $u$ to (\ref{1.2}) as a limit, and then show it
is complete.

Given $a\in(0,\infty)$, the existence of a unique solution 
$u_a\in C^\infty(\overline\Omega)$ 
to (\ref{1.2}) with $u_a = a$ on $\del \Ombar$ is easy
and well known; cf.~Chapter 14, \S{1} of \cite{T}.  The proof given there  
uses a simple combination of variational techniques and maximum principle
arguments.  Our strategy is to take $a\nearrow\infty$.  Thus we need to
consider how $u_a$ depends on $a$.

\begin{lemma} 
These solutions are monotonic in the parameter $a$:
\begin{equation}
a<b\Longrightarrow u_a\le u_b\ \mbox{ on }\ \Omega.
\label{2.1}
\end{equation}
\label{l2.1}
\end{lemma}
\noindent {\bf Proof:} Set $v=u_b-u_a$.  Then $v|_{\bW}=b-a>0$, while
\begin{equation}
\Delta v-\varphi_{ab}v=0,
\label{2.2}
\end{equation} 
with
\begin{equation}
\varphi_{ab}=\frac{e^{2u_b}-e^{2u_a}}{u_b-u_a}
=\frac{1}{u_b-u_a} \int_{u_a}^{u_b} 2 e^{2\sigma}\,d\sigma>0.
\label{2.3}
\end{equation}
Say $v_{{\mathrm min}}=v(p),\ p\in\Omega$.  If $v(p)<0$, then
$\Delta v(p)=\varphi_{ab}(p)v(p)<0$, which is impossible, so $v\ge 0$ 
on $\Omega$, yielding (\ref{2.1}).  \hfill $\Box$

$\text{}$ \newline
{\bf Remark 2.1.} This lemma and its variants will be used repeatedly in 
the following. In other useful variants, $a$ and $b$ can be 
functions on $\bW$, rather than just constants, and we can also
compare functions $u_a$ and $u_b$ that satisfy
\begin{equation}
\Delta u_a+k_a e^{2u_a}=K_0,\quad \Delta u_b+k_b e^{2u_b}=K_0.
\label{2.4}
\end{equation}
If $-k_a\ge -k_b>0$ on $\Ombar$, then (\ref{2.1}) holds.
\newline $\text{}$

The next step is to obtain a uniform upper bound for this
monotonic sequence.

\begin{lemma} There exists a locally bounded function $B$ on $\Om$ such that
\begin{equation}
e^{2u_a(p)}\le B(p),\quad \forall\ a<\infty.
\label{2.5}
\end{equation}
\label{l2.2}
\end{lemma}
\noindent {\bf Proof:} First consider the case where $\Om$ is a planar
domain, $\Om \subset \RR^2$. Define $\delta(p)={\mathrm{dist}}(p,\bW)$. 
We claim that 
\begin{equation}
e^{2u_a(p)}\le \frac{4}{\delta(p)^2},\quad \forall\ a<\infty.
\label{2.6}
\end{equation}
In fact, for any $\beta\in (0,\delta(p))$, let $D_\beta(p)$ be the disk 
of radius $\beta$, centered at $p$, with its Poincar\'e metric
\begin{equation}
g_{jk}=e^{2w}\delta_{jk}=\frac{4\beta^2}{(\beta^2-r^2)^2}\, \delta_{jk},\quad
r = |x-p|.
\label{2.7}
\end{equation}
Since $w$ solves (\ref{1.2}) and tends to $+\infty$ on $\del 
D_\beta(p)$, Lemma~\ref{l2.1} gives 
\begin{equation}
u_a\le w\ \mbox{ on }\ D_\beta(p).
\label{2.8}
\end{equation}
(\ref{2.6}) follows as $\beta\nearrow \delta(p)$. 

For the general case, use isothermal coordinates to get a neighbourhood 
$p\in \mathcal{O}_p \subset\Omega$ and a conformal map $\psi_p:D_1
\rightarrow \mathcal{O}_p$. We may assume that $\del\mathcal{O}_p$ is 
smooth and $\psi_p$ extends to a diffeomorphism on the closure. 
The Poincar{\'e} metric $g_p=e^{2w_p}g_0$ on $\mathcal{O}_p$ yields
a barrier, and hence an upper bound $w \leq B$, as above. It is clearly
possible to choose $B$ as a continuous function. \hfill $\Box$

\smallskip

Using Lemma \ref{l2.2}, we now let $a \nearrow \infty$ and obtain 
\begin{equation}
u_a(p)\nearrow u(p),\quad \forall\ p\in\Omega,\quad e^{2u(p)}\le
B(p).
\label{2.9}
\end{equation}
Each derivative of $u_a$ is locally uniformly bounded by elliptic regularity,
so convergence takes place in $C^\infty_{\mathrm{loc}}(\Omega)$, and hence
$u$ is a solution of (\ref{1.2}). 

$\text{}$ \newline
{\bf Remark 2.2.}  The upper bound (\ref{2.6}), valid for the limit $u$,
is not sharp as $p$ tends toward $\bW$.  An only slightly more involved
argument, using a conformal self-map of the disk $D_\beta(p)$, shows that
when $\bW$ is smooth, $\delta(p)^2 e^{2u(p)}\rightarrow 1$ as
$p\rightarrow \bW$.  A more refined result along these lines is the
content of Proposition \ref{p3.1} below.
\newline $\text{}$

It remains to demonstrate completeness.

\begin{lemma} Assume $\Ombar$ is a smoothly bounded and compact surface
with metric $g_0$, and $\Om$ is its interior. If $u$ is given as above, 
as the limit of the $u_a$, then the corresponding metric $g = e^{2u}g_0$
is complete on $\Om$. 
\label{l2.3}
\end{lemma}
\noindent {\bf Proof:} Let $\gamma:[0,L)\rightarrow\Omega$ be a unit-speed 
geodesic for $g$, with $L<\infty$, and suppose that $\gamma(t)$ does not
converge to a point in $\Omega$ as $t\rightarrow L$.  This curve also has 
finite length with respect to $g_0$, and so there exists $p\in \Ombar$ 
such that $\gamma(t)\rightarrow p$ as $t\to L$.  

As before, first consider the case where $\Omega$ is planar. Let 
$\mathcal{D}_p\subset\RR^2\setminus\Omega$ be a disk, tangent to $\bW$
at $p$. Regard $\RR^2$ as sitting inside the Riemann sphere $\widehat{\CC}$
and consider $\mathcal{D}'_p=\widehat{\CC}\setminus\mathcal{D}_p$, with its 
Poincar{\'e} metric $h$.  Thus $\Omega\subset\mathcal{D}'_p$.
The argument used in the first part of the proof of  Lemma~\ref{l2.2},
applied to a sequence of disks decreasing to $\mathcal{D}'_p$, also gives
\begin{equation}
g\ge h\ \mbox{ on }\ \Omega.
\label{2.10}
\end{equation}
But $(\mathcal{D}'_p,h)$ is complete, so the $h$-length of 
$\gamma$ is infinite, and hence $\Omega$ is complete with respect to $g$.

To handle the general case, assume that $\Ombar$ sits inside
a slightly larger (open) Riemann surface $M$ and the metric $g_0$
is extended smoothly. If $\gamma(t) \to p \in \bW$ as $t \nearrow L$,
choose a small holomorphic disk $\mathcal{D}$ containing $p$, roughly cut in 
half by $\bW$. Choose a smooth curve in $\mathcal{D} \setminus\Ombar$ 
hitting $\bW$ transversally at $p$ and let $p_j\rightarrow p$ along this 
curve.  Denote by $e^{2v_j}g_0$ the Poincar{\'e} metric on $\mathcal{D}
\setminus\{p_j\}$ obtained by pulling back (\ref{1.5}), and let  
$\mathcal{O}\subset\subset\mathcal{D}$ be a smaller disk containing $p$
and the sequence $p_j$. We see that, for each $j<\infty$, there exists 
$A(j)<\infty$ such that
\begin{equation}
u_a\ge v_j\ \mbox{ on }\ \bW\cap\mathcal{O},\ \mbox{ for }\ a\ge A(j).
\label{2.11}
\end{equation}
Also, considering $u_1$ (which equals $1$ on $\bW$) we see that there exists
$B > 0$ such that $u_1\ge v_j-B$ on $\pa\mathcal{O}\cap\Omega$, 
for all $j$, hence
\begin{equation}
u_a\ge v_j-B\ \mbox{ on }\ \pa\mathcal{O}\cap\Omega,\ \mbox{ for }\ a\ge 1.
\label{2.12}
\end{equation}
Hence, by Remark 2.1,
\begin{equation}
u_a\ge v_j-B\ \mbox{ on }\ \mathcal{O}\cap\Omega,\ \mbox{ for }\ a\ge
\max(1,A(j)),
\label{2.13}
\end{equation}
so
\begin{equation}
u\ge v_j-B\ \mbox{ on }\ \mathcal{O}\cap\Omega,\quad \forall\ j.
\label{2.14}
\end{equation}
Hence
\begin{equation}
u\ge v-B \ \mbox{ on }\ \mathcal{O}\cap\Omega,
\label{2.15}
\end{equation}
where $e^{2v}\,g_0$ is (\ref{1.5}) pulled back to $\mathcal{D}
\setminus\{p\}$. This is enough to give completeness.
\hfill $\Box$

Putting these lemmas together we obtain
\begin{proposition} If $\Om$ is the interior of a smooth two-dimensional 
manifold with boundary 
$\Ombar$, then $\Om$ admits a Poincar\'e metric. 
\label{p2.4}
\end{proposition} 

As we have emphasized, Proposition \ref{p2.4} implies $\Omega$ is 
holomorphically covered by $D_1$.  In particular, if $\Omega$ is simply
connected then there exists a holomorphic diffeomorphism $\varphi:D_1
\rightarrow \Omega$.  It is useful to recall the linear PDE treatment 
of this result, in which one picks $p\in\Omega$ and takes the Green function
$u\in H^{1-\ep}(\Omega)\cap C^\infty(\Ombar\setminus\{p\})$, satisfying
\begin{equation}
\Delta u=2\pi \delta_p,\ \text{ on }\ \Omega,\quad
u\bigr|_{\bW}=0.
\label{2.16}
\end{equation}
Then $u(x)$ behaves like $\log |x|$ in local normal coordinates centered
at $p$, and the hypothesis that $\Omega$ is simply connected implies that
there exists a harmonic conjugate $v$, smooth and well defined mod
$2\pi\ZZ$, on $\Ombar\setminus\{p\}$, and the function
\begin{equation}
\Phi(x)=e^{u(x)+iv(x)}
\label{2.17}
\end{equation}
yields a holomorphic map $\Phi:\Omega\rightarrow D_1$, extending to a 
smooth map $\Phi:\Ombar\rightarrow \overline{D}_1$.  It follows from the
Hopf lemma (also known as Zaremba's principle) that $\pa_\nu u<0$ on 
$\bW$, which via the Cauchy-Riemann equations implies $\Phi$ maps $\bW$
locally diffeomorphically onto $S^1=\pa D_1$.  Now the argument principle 
implies the degree of the map $\bW\rightarrow S^1$ equals the number of
preimages of any $q\in D_1$ (counting multiplicity).  This number is
clearly one for $q=0$, so $\Phi$ is the desired holomorphic
diffeomorphism.  This argument gives us something extra; $\Phi$ extends to
a smooth diffeomorphism of $\Ombar$ onto $\overline{D}_1$.  This has 
implications for the boundary behavior of the Poincar{\'e} metric on
$\Omega$, which we will explore and extend in the next section.

\section{Boundary regularity when $\Ombar$ is smooth}

In this section we analyze the boundary behavior of the function
$u\in C^\infty(\Omega)$ providing the Poincar{\'e} metric $e^{2u}g_0$,
when $\Omega$ is smoothly bounded.  Throughout this section 
we let $x$ denote the distance function (with respect to $g_0$), which
is well-defined and smooth in some sufficiently small neighbourhood of 
the boundary $\bW$, and shall often also use $y$ as a local coordinate 
along $\bW$. 

\begin{proposition} Suppose that $\Ombar$ is smoothly bounded
and compact, with Poincar{\'e} metric $e^{2u}g_0$.
Then as $x \to 0,\ u$ has an asymptotic expansion of the form
\begin{equation}
u(x,y) \sim \log (1/x) + u_1(y) x + u_2(y)x^2 + \cdots,
\label{eq:asexpu}
\end{equation}
where the coefficient functions $u_j(y)$ all lie in $\calC^\infty(\bW)$.
Equivalently,
\begin{equation}
e^{-u}\in C^\infty(\Ombar),\ \text{ and }\ \pa_\nu e^{-u} \equiv 1.
\label{3.2}
\end{equation}
\label{p3.1}
\end{proposition}

This is the direct analogue of the expansion valid for solutions 
of the singular Yamabe problem in higher dimensions (at least
in the most favourable case), cf.\ \cite{Ma2}. The proof has 
two steps: first, barrier techniques are used to obtain rough
(scale-invariant) estimates for the solution, and at that
point some techniques from the linear analysis of \cite{Ma1} are
used to improve this to full tangential regularity and an expansion.

$\text{}$ \newline
{\bf Remark 3.1.} We note that the second condition in (\ref{3.2})
is an automatic consequence of the first.  In fact, $W=e^{-u}$ satisfies
\begin{equation}
\Delta W=\frac{|\nabla W|^2-1}{W}-K_0W,\quad W\bigr|_{\bW}=0.
\label{3.3}
\end{equation}
If $W\in C^\infty(\Ombar)$, then the right side of (\ref{3.3}) must be
continuous on $\Ombar$, which implies $\pa_\nu W|_{\bW}=1$.
\newline $\text{}$

$\text{}$ \newline
{\bf Remark 3.2.} In case $\Ombar$ is smoothly bounded and simply
connected, the smoothness of $e^{-u}$ on $\Ombar$ is a simple consequence
of the fact that the holomorphic diffeomorphism $\Phi:\Omega\rightarrow
D_1$ given by (\ref{2.17}) extends to a smooth diffeomorphism $\Phi:
\Ombar\rightarrow \overline{D}_1$, plus the fact that
$e^{2u}g_0=\Phi^*(g_P)$, where $g_P$ is the Poincar{\'e} metric on $D_1$.
\newline $\text{}$

$\text{}$ \newline
{\bf Remark 3.3.} The smoothness condition in (\ref{3.2}) is clearly
invariant when $g_0$ is replaced by $g_1=e^{2w}g_0$ with $w\in C^\infty
(\Ombar)$, and hence so are the rest of the conclusions in Proposition 
\ref{p3.1}.
\newline $\text{}$

We will implement Remark 3.3 using the following result.

\begin{lemma} For each connected component $\gamma$ of $\bW$, there is
a collar neighborhood ${\cal{C}}$ and a $C^\infty$ conformal
diffeomorphism $\varphi:\overline{\cal{C}}\rightarrow \overline{A}_b$
onto an annulus $\overline{A}_b=\{z\in\CC:e^{-\pi/b}\le |z|\le 1\}$,
as in (\ref{1.7}).
\label{l3.2}
\end{lemma}

\noindent
{\bf Proof:} Taking a collar neighborhood ${\cal{C}}_0$ of $\gamma$,
we can produce a simply connected ${\cal{O}}$ with smooth boundary 
such that a collar neighborhood of $\pa{\cal{O}}$ is identified with
${\cal{C}}_0$.  Then we can apply the construction mentioned at the
end of \S{2}, obtaining a $C^\infty$ conformal diffeomorphism
$\Phi:\overline{\cal{O}}\rightarrow \overline{D}_1$.  The inverse image
of $\overline{A}_b$, for $b$ sufficiently large, can then be identified
with the desired collar neighborhood $\overline{\cal{C}}$ of $\gamma$.
\hfill $\Box$
\newline $\text{}$

Using Lemma \ref{l3.2}, we can construct $w\in C^\infty(\Ombar)$ such that 
$g_1=e^{2w}g_0$ has the property that each boundary component $\gamma$
of $\bW$ has a collar neighborhood that is {\it isometric} to
$\overline{A}_b$.  We now renotate, giving $g_1$ the label $g_0$.
We are ready to establish the following estimate.

\begin{lemma} In the setting of Proposition \ref{p3.1}, we have, near
$\bW$,
\begin{equation}
u=\log \frac{1}{x}+v,\quad |v|\le Cx.
\label{3.4}
\end{equation}
\label{l3.3}
\end{lemma}

\noindent
{\bf Proof:} Recall the conformal diffeomorphism
$\varphi:\overline{\cal{C}}\rightarrow \overline{A}_b$ constructed in 
Lemma \ref{l3.2}.  Pulling the Poincar{\'e} metric (\ref{1.7}) on $A_b$
back via $\varphi$ produces the Poincar{\'e} metric on ${\cal{C}}$, say
$e^{2u_1}g_0$, and we clearly have
\begin{equation}
u\le u_1,\ \text{ on }\ {\cal{C}}.
\label{3.5}
\end{equation}
It is clear from the representation (\ref{1.7}) that this produces an
upper bound on $u$ of the form asserted in (\ref{3.4}).

It remains to produce an appropriate lower bound on $u|_{\cal{C}}$.
It is equivalent to produce a lower bound on the metric
$e^{2u}g_0|_{\cal{C}}$, pulled back to $A_b$ via $\varphi$.  To accomplish
this, we supplement the family of metrics (\ref{1.7}) with the following
family of metrics on $A_b$:
\begin{equation}
g^\beta_{jk}=\Bigl[ \frac{\beta}{r\,\text{sinh}(\beta\,\log\frac{1}{r})}
\Bigr]^2 \delta_{jk}.
\label{3.6}
\end{equation}
These are obtained by analytically continuing (\ref{1.7}) to purely
imaginary $b$.  A direct check shows that these metrics satisfy 
(\ref{1.2}) on $\{x\in\RR^2:0<|x|<1\}$, for each $\beta\in (0,\infty)$,
and the metric (\ref{1.5}) arises as the limit as $\beta \searrow 0$.
The metrics (\ref{3.6}) are complete at the outer boundary $\{z:|z|=1\}$,
and, given any fixed $b>0$, for large $\beta$ they are quite small on the
inner boundary $\{z:|z|=b\}$ of $A_b$.  Choosing $\beta$ sufficiently 
small gives the desired lower bound, establishing (\ref{3.4}).
\newline $\text{}$

What we have accomplished thus far is to show that the conformal factor 
$u$ giving the Poincar\'e metric $g = e^{2u}g_0$ may be written 
on a collar neighborhood ${\cal{C}}$ of each boundary component as 
$u = \log(1/x) + v$, where $|v| \leq Cx$ for $0 < x \leq x_0$. Here 
$x = 1-r$, $r = |z|$ on the annulus $A_b$, identified with ${\cal{C}}$.
Notice that $u$ satisfies
\[
\Delta_{g_0} u - e^{2u} = 0
\]
in $A_b$, and in addition, letting $y$ be the polar angular variable 
$\theta$ on $A$, then $\Delta_{g_0} = \del_x^2 - (1-x)^{-1}\del_x + 
(1-x)^{-2}\del_y^2$ there. Hence $\Delta_{g_0}(\log(1/x)) - e^{2\log(1/x)} 
= 1/x(1-x)$, and so $\Delta_{g_0}v + (1/x^2)(1-e^{2v}) = 1/x(1-x)$; we
rewrite this finally as
\begin{equation}
Lv \equiv (x^2 \Delta_{g_0} - 2)v = Q(v) + r(x)
\end{equation}
where $r(x) = - x/(1-x)$ and $Q(v) = e^{2v} - 1 - 2v$ is smooth and 
vanishes quadratically as $v \to 0$.

It may seem that we have lost ground since the linear operator
$L = x^2 \Delta_{g_0} - 2$ appearing here, while elliptic in the interior, 
is uniformly degenerate at $\bW$. However, this sort of degenerate
elliptic operator is well-understood, and \cite{Ma1} contains 
a general framework for studying degenerate operators of this type.
We now state the results we need from that paper and then apply
them to our purposes.

We shall use a scale of weighted H\"older spaces,
$x^\gamma\Lambda^{\ell,\alpha,\ell'}_0(\Ombar)$ for $\gamma \in \RR$ 
and $\ell, \ell' \in {\mathbb N}$, $\ell \geq \ell'$. First, when 
$\gamma=0$ and $\ell' = 0$, then $\Lambda^{\ell,\alpha,0}_0
\equiv \Lambda^{\ell,\alpha}_0$ is the `geometric H\"older space' 
associated to the covariant derivative for the metric $g_1 = x^{-2}g_0$ 
(or any metric smoothly quasi-isometric to this). This mean that $w$
is in this space if the supremum over all $g_1$-unit balls 
of the H\"older seminorms with exponent $\alpha$ of the functions 
$(x\del_x)^j (x\del_y)^k w$, $j+k \leq \ell$, 
is finite; the norm is the obvious one. Note that all derivatives
here are taken with respect to the degenerate vector fields $x\del_x$ 
and $x\del_y$. The space $\Lambda^{\ell,\alpha,\ell'}_0$, still with 
weight parameter $\gamma = 0$, consists of those elements $w \in 
\Lambda^{\ell,\alpha}_0$ such that $\del_y^k w \in
\Lambda^{\ell-k,\alpha}_0$ for $0 \leq k \leq \ell'$. In other words,
up to $\ell'$ of the $x\del_y$ derivatives may be replaced by
derivatives with respect to the nondegenerate vector field $\del_y$. Finally,
\begin{equation}
x^\gamma\Lambda^{\ell,\alpha,\ell'}_0 = \{w = x^\gamma \tilde{w}:
\tilde{w} \in \Lambda^{\ell,\alpha,\ell'}_0\}.
\end{equation}
Clearly
\begin{equation}
L: x^\gamma \Lambda^{\ell+2,\alpha,\ell'}_0 \to x^\gamma \Lambda^{\ell,
\alpha,\ell'}_0
\label{eq:mpgs}
\end{equation}
is bounded for every $\gamma \in \RR$ and $0 \leq \ell' \leq 
\ell$. But this map may be ill-behaved in various ways, and to understand
this we must compute the indicial roots of $L$. By definition,
$\gamma$ is an indicial root of $L$ if $L(x^\gamma) = O(x^{\gamma+1})$. 
But 
\begin{equation}
L(x^\gamma) = x^2 \Delta_{g_0}x^\gamma - 2 x^\gamma =
(\gamma^2 - \gamma - 2)x^\gamma + O(x^{\gamma + 1}),
\end{equation}
so this can only happen if $\gamma = \gamma_\pm$, where $\gamma_- = -1$
and $\gamma_+ = 2$. These are the only two indicial roots of $L$.

It is not hard to check that (\ref{eq:mpgs}) fails to have closed range 
when $\gamma = \gamma_\pm$. On the other hand, Corollary 6.4 and 
Proposition 5.30 in \cite{Ma1} give 
\begin{lemma} The map (\ref{eq:mpgs}) is Fredholm
of index zero when $-1 < \gamma < 2$. 
\end{lemma}

The proof relies on the construction of a parametrix $G$ for $L$
such that
\begin{equation}
G: x^\gamma\Lambda^{\ell,\alpha,\ell'}_0 \longrightarrow 
x^\gamma\Lambda^{\ell+2,\alpha,\ell'}_0
\label{3.10}
\end{equation}
is bounded for all $0 \leq \ell' \leq \ell$, and such that both $GL - I$ 
and $LG - I$ are compact. This uses the restriction $-1<\gamma < 2$, 
and immediately implies that (\ref{eq:mpgs}) is Fredholm when
$\gamma$ is in this range. 
The vanishing of the index follows from the formal self-adjointness of 
$L$ (or alternately, because $L$ is real and scalar). 

To proceed further, we also need a regularity theorem, which is 
Proposition 3.28 in \cite{Ma1}:
\begin{lemma} Suppose $-1 < \gamma \leq N$ and $Lw = f$,
where $w \in x^\gamma L^\infty$ and $f = x^N \tilde{f}$,
$\tilde{f} \in C^\infty(\Ombar)$, then necessarily
$w = x^2 w_1 + x^N w_2 + x^2 \log x \, w_3$ where 
$w_1, w_2, w_3 \in C^\infty(\Ombar)$. If $N > 2$ here
then $w_3 = 0$, i.e., the expansion for $u$ has no logarithms. 
\label{l3.4}
\end{lemma}
\noindent {\bf Remark 3.4.} This is a specialization of a more general 
result which, for this operator $L$, states that if $w \in x^\gamma 
L^\infty$, $\gamma > -1$, and $Lw = f$ where $f$ has a general 
polyhomogeneous expansion 
(with all exponents greater than or equal to $\gamma$), then $w$
also has a polyhomogeneous expansion of the same form, although possibly
with terms with extra logarithmic factors.
\newline $\text{}$

This lemma applies immediately as follows: if $w\in x^\gamma\Lambda^{\ell,
\alpha, \ell'}_0$ with $-1<\gamma<2$, and $Lw = 0$ then in particular, $w$ 
is smooth on $\Ombar$ and vanishes at
$\bW$. Since solutions of $Lw=0$ satisfy the maximum principle, we get
$w=0$. Hence (\ref{eq:mpgs}) is injective, and thus an isomorphism, when
$-1 < \gamma < 2$. 

Now recall the decomposition $u = \log(1/x) + v$ where $|v| \leq Cx$, 
i.e., $v \in x L^\infty$. Since $v$ satisfies a semilinear elliptic 
equation which is uniformly elliptic in unit balls relative to the
metric $g_1 = x^{-2}g_0$, we may use standard Schauder estimates in 
each of these balls, and recall the initial definition of the weighted
H\"older spaces with $\ell'=0$ to conclude that $v \in x 
\Lambda^{\ell,\alpha}_0$ for every $\ell \geq 0$.  

Our next (and final) major claim is that $v \in x\Lambda^{\ell,\alpha,
\ell'}_0$ for every $0 \leq \ell' \leq \ell$. Set $\gamma = 1$ and let
$G$ denote the corresponding inverse for $L$. Write the equation
$Lv = Q(v) + r(x)$ as $v = GQ(v) + G(r)$; this is legitimate
because both $Q(v)$ and $r(x)$ lie in $x\Lambda^{\ell,\alpha}_0$ for
every $\ell \geq 0$. In fact, since $r(x) \in x\Lambda^{\ell,\alpha,
\ell'}_0$ for every $\ell \geq \ell'$, the boundedness of (\ref{3.10}) 
shows that the final term is completely tangentially regular. 
Next, write $v = x \tv$, so that 
$\tv \in \Lambda^{\ell,\alpha}_0$. Then $Q(v) = x^2 \tQ(x,\tv)$,
where $\tQ(x,s)$ again vanishes quadratically as $s \to 0$. 
Let us make the inductive hypothesis that $v \in x\Lambda^{\ell,\alpha,
\ell'}_0$ for some fixed $\ell'$ (and every $\ell \geq \ell'$). This
is clearly true when $\ell' = 0$, so we must show that if it
true for some value of $\ell'$, then it is true when $\ell'$
is replaced by $\ell'+1$. This uses a commutator argument. 
In fact, 
\[
\del_y v = \del_y GQ(v) + \del_yG(r).
\]
Neglecting the final term on the right, which we already know
has the correct regularity, reexpress the other term on the right as
\[
G(\del_y x^2 \tQ(\tv)) + [\del_y,G]Q(v). 
\]
By Proposition~3.30 in \cite{Ma1}, the commutator $[\del_y,G]$ 
enjoys the same mapping properties (\ref{3.10}) as $G$ itself, 
and so the second term here lies in $x\Lambda^{\ell,\alpha,\ell'}$, 
by the inductive hypothesis. On the other hand, $\del_y(x^2\tQ(\tv))
= (x\del_y)x\tQ(\tv)$, and since $x\tQ(\tv) \in x\Lambda^{\ell,
\alpha,\ell'}_0$, we see that this first term also has this 
same regularity. Thus all terms in this expression for $\del_y v$ 
lie in $x\Lambda^{\ell,\alpha,\ell'}_0$, and so $v \in x \Lambda^{\ell,
\alpha,\ell'+1}_0$ for all $\ell \geq \ell'+1$. 
This proves the claim, and shows that $v$ is fully tangentially regular.

It remains to establish that $v= x\tilde{v}$ where $\tilde{v} \in 
C^\infty(\Ombar)$. One extra consideration we need 
to address is that there are no logarithmic terms in the
expansion for $v$, the presence of which might be suspected 
from Lemma~\ref{l3.4}.
Define $\calA^\nu$ to be the intersection of $x^\nu
\Lambda^{\ell,\alpha,\ell'}_0$ over all $0 \leq \ell' \leq \ell < \infty$; 
we have shown that $v \in \calA^1$. To deduce its expansion, write 
$L = x^2 \del_x^2 - 2 + E$, where $E$ consists of all `error terms' 
(i.e., $x^2 \del_y^2$ and $-x^2/(1-x)\del_x$). Now regard the 
equation for $v$ as an ODE in $x$ with values in functions 
smooth on the boundary:
\begin{equation}
x^2 \del_x^2 v - 2v = -Ev + Q(v) + r(x).
\end{equation}
We think of the whole right hand side as an inhomogeneous term.
Recall that $r(x) = -x/(1-x) = -x - x^2 - \ldots$, and
$Q(v) = e^{2v} - 1-2v = 2v^2 + O(v^3)$. Then at the 
first stage the right hand side has the form $-x + f_2$, with
$f_2 \in \calA^2$. Integrating the ODE gives $v = -(1/2)x + 
v_2$, $v_2 \in \calA^2$. Inserting this into the right side
shows that the sum of these terms on the right have the
form $-x + f_3$, $f_3 \in \calA^3$. The fact that the $x^2$ term
in this expansion vanishes is a special feature, due to a
fortuitous cancellation; the absence of this term is what
precludes the logarithm terms in the expansion for $v$.
Integrating the ODE again shows that $v = -(1/2)x + v_2(y) x^2
+ v_3$, where $v_3 \in \calA^3$ and $v_2(y) \in C^\infty(\bW)$.
Inserting this back into the right side and iterating this 
argument gives the complete expansion for $v$. This completes the proof. 
\hfill $\Box$

$\text{}$ \newline
{\bf Remark 3.5.} From (\ref{3.3}) one can compute $\pa^2_\nu W|_{\bW}$,
and see that in the expansion (\ref{eq:asexpu}) $u_1(y)=\kappa(y)/2$, where
$\kappa(y)$ is the curvature of $\bW$ at $y$.  On the other hand,
the coefficient $u_2(y)$ depends on the global behavior of $\Omega$,
as one can see by examining (\ref{1.7}) for different values of $b$.

\section{General planar domains}

In this section we construct Poincar{\'e} metrics on general planar
domains, as long as the complement contains at least two points.
To begin, given $\Omega\subset\RR^2$, open and connected, take 
a sequence $\Omega_\nu$ bounded, 
with smooth boundary, such that $\Omega_\nu\subset\subset
\Omega_{\nu+1}$ and $\Omega_\nu\nearrow\Omega$, in the sense that any
compact $K\subset\Omega$ is contained in $\Omega_\nu$ for large $\nu$.
Let $u_\nu$ be the solutions to (\ref{1.2}) on $\Omega_\nu$ such that 
$u_\nu|_{\bW_\nu}=+\infty$ and $g^\nu_{jk}=e^{2u_\nu}\delta_{jk}$ are
complete metric tensors on $\Omega_\nu$, as in Proposition \ref{p2.4}
(obtained as in (\ref{2.9})).  The argument
used to prove Lemma \ref{l2.1} shows that $u_\nu\searrow$ as $\nu\nearrow\infty$.
Our main goal in this section is to establish the following.

\begin{proposition}  If $\Omega\subset\RR^2$ is a connected open set 
with the property that $\RR^2\setminus\Omega$ contains at least two points,
then 
\begin{equation}
u_\nu\searrow u\ \mbox{ as }\ \nu\nearrow\infty,
\label{4.1}
\end{equation}
where $u\in C^\infty(\Omega)$, solving (\ref{1.2}), and the metric tensor 
$g_{jk}=e^{2u}\delta_{jk}$ is a complete metric tensor on $\Omega$, of Gauss 
curvature $-1$.
\label{p4.1}
\end{proposition}

As a warm-up, we first give a simple proof of the following special case,
which extends Proposition \ref{p2.4}, in the case of planar domains.

\begin{proposition} Proposition \ref{p4.1} holds when $\Omega\subset\RR^2$
is a bounded, connected, open set, whose boundary satisfies the following
regularity hypothesis:
\begin{equation}
\mbox{Each }p\in\bW\mbox{ is the endpoint of a line segment in }\RR^2
\setminus\Omega.
\label{4.2}
\end{equation}
\label{p4.2}
\end{proposition}
\noindent {\bf Proof:} 
First we need to get a bound on $u_\nu$ from below.  Indeed,
taking $\Omega$ inside a sufficiently large disk $D_b(0)$, with 
Poincar{\'e} metric $e^{2w}\delta_{jk}$, then $u_\nu\ge w$ on $\Omega_\nu$.  
This gives a locally bounded $u$ on $\Omega$ which satisfies (\ref{4.1}).
As before elliptic estimates give smooth convergence to $u\in 
C^\infty(\Omega)$, solving (\ref{1.2}).  

Completeness remains to be demonstrated.  Under the hypothesis (\ref{4.2}), 
the completeness proof goes as follows.  Say $\gamma:
[0,L)\rightarrow\Omega$ is a unit-speed geodesic (for $g_{jk}$) and suppose
$L<\infty$ and $\gamma(t)$ does not converge to a point in $\Omega$ as
$t\rightarrow L$.  As in the proof of Proposition 2.3, we have $\gamma(t)
\rightarrow p$ for some $p\in\bW$.  

Let $\ell$ be a line segment in $\RR^2\setminus\Omega$ with $p$ as an 
endpoint.  Regard $\RR^2\subset\widehat{\CC}$.  
Now it is elementary to produce a conformal diffeomorphism 
$\psi:\widehat{\CC}\setminus\ell\rightarrow D_1$; 
pull back the Poincar{\'e} metric on $D_1$ to get a complete metric 
$e^{2w}\delta_{jk}$ on $\widehat{\CC}\setminus\ell$, of Gauss curvature $-1$.
Again we have $u_\nu\ge w$ on $\Omega_\nu$, for each $\nu<\infty$, and hence
$u\ge w$ on $\Omega$, and the completeness of 
$e^{2u}\delta_{jk}$ on $\Omega$ is proven.
\hfill $\Box$

The proof of Proposition \ref{p4.1} for more general $\Omega$ requires more 
work, which we now undertake. To get a lower bound on $u_\nu$ this time, 
we make use of the following result.

\begin{lemma} The region $\CC\setminus\{0,1\}$ has a Poincar{\'e} metric.
\label{l4.3}
\end{lemma}

We will give a curvature equation
proof of this lemma after we apply it to prove Proposition \ref{p4.1}.

Returning to the estimation of $u_\nu$ in the proof of Proposition 
\ref{p4.1}, say $p_1,p_2\in\RR^2\setminus\Omega$.  Lemma \ref{l4.3} also holds 
for $\RR^2\setminus\{p_1,p_2\}$, which therefore has a 
Poincar\'e metric
\begin{equation}
h_{jk}=e^{2w}\,\delta_{jk}.
\label{4.3}
\end{equation}
Now as in Lemma \ref{l2.1} we have $u_\nu\ge w$ on 
$\Omega_\nu$.  Hence, as before we can deduce that $u_\nu\searrow u$
with $u\in C^\infty(\Omega)$ satisfying (1.1).  To prove that the metric
$e^{2u}\delta_{jk}$ is complete on $\Omega$, we argue as before that if
not, there would exist a unit-speed $\gamma:[0,L)\rightarrow\Omega$ with
$\gamma(t)\rightarrow p\in\bW$ as $t\rightarrow L$.  Here $\bW$ denotes
the boundary of $\Omega$ in $\widehat{\CC}$, so either $p\in\RR^2\setminus
\Omega$ or $p=\infty$.  We now bring in the metric (\ref{4.3}), with $p_j\in
\RR^2\setminus\Omega$ and with $p_1=p$ if $p\neq\infty$; using $u\ge w$
we again have that $e^{2u}\delta_{jk}$ is complete on $\Omega$.
This gives Proposition \ref{p4.1}, modulo a proof of Lemma \ref{l4.3}.

We turn now to a proof of Lemma \ref{l4.3}.  One ingredient will be a metric 
on $\CC\setminus\{0,1\}$ of the form $e^{2w_0(z)}\delta_{jk}$, with 
\begin{equation}
e^{w_0}=A\, \frac{(1+r^a)^b}{r^c}\, \frac{(1+\rho^a)^b}{\rho^c},\quad
r=|z|,\ \rho=|z-1|,
\label{4.4}
\end{equation}
with $A,a,b,c>0$.  A calculation of the Gauss curvature for this metric 
gives
\begin{equation}
K=- \frac{a^2b}{A^2} \Bigl[
\frac{r^{a-2+2c}\rho^{2c}}{(1+r^a)^{2+2b}(1+\rho^a)^{2b}}+
\frac{r^{2c}\rho^{a-2+2c}}{(1+r^a)^{2b}(1+\rho^a)^{2+2b}}\Bigr].
\label{4.5}
\end{equation}
We have $K<0$ and it is bounded away from zero as long as
\begin{equation}
a-2+2c\le 0,\quad 4c-a-4ab-2\ge 0.
\label{4.6}
\end{equation}
For example, we can take
\begin{equation}
a=\frac{1}{3},\quad b=\frac{1}{2},\quad c=\frac{5}{6},
\label{4.7}
\end{equation}
the parameters used in \cite{Kr}, pp.~78--80.  If $A>0$ is small enough,
we have $K\le -1$.

Fix such $A$, let $\Omega_\nu\nearrow\Omega=\CC\setminus\{0,1\}$, and
take $u_\nu\in C^\infty(\Omega_\nu)$, as in Proposition \ref{p4.1}, with $u_\nu
\searrow$ as $\nu\nearrow\infty$.  A variant of the proof of Lemma \ref{l2.1}
gives $u_\nu\ge w_0$ on $\Omega_\nu$, with $w_0$ given by (\ref{4.4}), so
we have convergence: $u_\nu\rightarrow u$ with $u\in C^\infty(\Omega)$
satisfying (\ref{1.2}) and $u\ge w_0$ on $\Omega$.  However, the metric
$e^{2w_0}\delta_{jk}$ is not complete, so we need to do some more work
to show that $e^{2u}\delta_{jk}$ is complete on $\CC\setminus\{0,1\}$.

To check completeness of $e^{2u}\delta_{jk}$ near $0$, we compare it with 
the metric (\ref{1.5}), i.e., $e^{2v}\delta_{jk}$, where
\begin{equation}
e^{2v}=\Bigl(r\, \log \frac{1}{r}\Bigr)^{-2},
\label{4.8}
\end{equation}
on $0<r<1$.  Given the convergence $u_\nu\rightarrow u$, we can find a 
constant $B\ge 0$ such that $u_\nu\ge v-B$ on $\{z\in\CC:|z|=1/2\}$.
Note that $e^{2(v-B)}\delta_{jk}$ has curvature $-e^{2B}\le -1$ on $D^*$.
Now a variant of Lemma \ref{l2.1} gives $u_\nu\ge v-B$ on $\{z\in\Omega_\nu:
|z|\le 1/2\}$, and hence
\begin{equation}
u\ge v-B\ \mbox{ on }\ \{z:0<|z|\le 1/2\}.
\label{4.9}
\end{equation}
This implies completeness of $e^{2u}\delta_{jk}$ near $0$.  Completeness
near $1$ is established similarly.  The formula (\ref{1.5}) also defines a
Poincar{\'e} metric on $\{z:|z|>1\}$, and this can be used to show that
$e^{2u}\delta_{jk}$ is complete near $\infty$.  Lemma \ref{l4.3} is proven.

$\text{}$ \newline
{\bf Remark.}  Lemma \ref{l4.3} is equivalent to the assertion that there is a 
holomorphic covering map
\begin{equation}
\psi:D_1\longrightarrow \CC\setminus\{0,1\}.
\label{4.10}
\end{equation}
This result is an ingredient in the classical theorems of Picard.  
The map $\psi$ can be constructed explicitly via elliptic function 
theory.  Cf.~Chapter 7 of \cite{Ahl}; this provided the original proof.
This covering can also be constructed by applying Schwarz reflection to the 
Riemann mapping function of a special domain on $\CC$ (cf.~Chapter 5, \S{6}
of \cite{T}).  A variant of (\ref{4.4}), obtained by adding multiples of 
(\ref{1.5}) and its images near $0,\ 1$, and $\infty$, was produced in 
\cite{GR} and shown there to have Gauss curvature $\le -1$ and to be 
complete on $\CC\setminus\{0,1\}$; cf.~\cite{Kob}, pp.~7--10.
\newline $\text{}$

We can produce other Riemann surfaces covered by the disk, using the 
following simple result.
\begin{proposition} If $M$ is a Riemann surface with a holomorphic
covering map $\psi:D_1\rightarrow M$ and $\Omega\subset M$ is a nonempty
open connected set, then there exists a holomorphic covering map
$\varphi:D_1\rightarrow \Omega$.
\end{proposition}
\noindent {\bf Proof:} If $\mathcal{O}\subset D_1$ is a connected component of $\psi^{-1}
(\Omega)$, then $\psi$ restricts to a holomorphic covering $\psi:\mathcal{O}
\rightarrow\Omega$.  By Proposition \ref{p4.2}, there exists a holomorphic covering
$\tilde{\psi}:D_1\rightarrow\mathcal{O}$.  Composing gives the holomorphic
covering $\varphi:D_1\rightarrow\Omega$.
\hfill $\Box$

We will not dwell on applications of this last proposition, since they
would all be subsumed by the results of \S{6}.

\section{Koebe's disk theorem}

Here we make note of a simple curvature proof of some
results of P.~Koebe on the family $\cal{S}$ of univalent (i.e., one-to-one)
holomorphic maps $f:D_1\rightarrow\CC$ satisfying $f(0)=0,\ f'(0)=1$.
Here is the first result.

\begin{proposition} There exists a constant $b\in (1,\infty)$ such that
for any $f\in{\cal{S}},\ \Omega=f(D_1)$ has the property
\begin{equation}
\frac{1}{b}\le \text{dist}(0,\bW)\le 1.
\label{5.1}
\end{equation}
\label{p5.1}
\end{proposition}

Here (and in (\ref{5.2}), (\ref{5.4}) below) we use Euclidean distance, so
that $\text{dist}(0,\bW)=\inf \{|z|:z\in\bW\}$.  L.~Bieberbach showed that 
one can take $b=4$, and this is sharp.  This sharpened version of 
Proposition \ref{p5.1} is called the Koebe-Bieberbach quarter theorem.  
Our method does not yield $b=4$.  The following result is equivalent to
Proposition \ref{p5.1}.

\begin{proposition} Let $\Omega$ be a proper, simply connected domain 
in $\CC$.  Let $e^{2u}\, \delta_{jk}$ be the Poincar{\'e} metric on 
$\Omega$.  Then, for all $p\in\Omega$,
\begin{equation}
\frac{1}{2}\ \text{dist}(p,\bW)\le e^{-u(p)}\le \frac{b}{2}\,
\text{dist}(p,\bW).
\label{5.2}
\end{equation}
\label{p5.2}
\end{proposition}

To see the equivalence, note that if $f:D_1\rightarrow \Omega$ is
biholomorphic and $\gamma(z)^2 |dz|^2$ is a metric on $\Omega$, 
then $D_1$ inherits the metric $\gamma(f(z))^2 |f'(z)|^2 \,|dz|^2$. 
Thus the Poincar{\'e} metric $e^{2u}\,\delta_{jk}$ induced on $\Omega$ has 
\begin{equation}
e^{-u(f(z))}=\frac{1}{2} (1-|z|^2) |f'(z)|.
\label{5.3}
\end{equation}
Picking a biholomorphic $f$ such that $f(0)=p$ yields the equivalence
of these propositions easily.  In addition, comparing (\ref{5.2}) 
and (\ref{5.3}) gives the following result.

\begin{proposition} If $f:D_1\rightarrow \Omega$ is a biholomorphic
map, then, for all $z\in D_1$,
\begin{equation}
\text{dist}(f(z),\bW)\le (1-|z|^2)|f'(z)|\le b\, \text{dist}(f(z),\bW).
\label{5.4}
\end{equation}
\label{p5.3}
\end{proposition}

We note that the upper estimate of $\text{dist}(0,\bW)$ in (\ref{5.1}) and
(equivalently) the lower estimate on $e^{-u(p)}$ in (\ref{5.2}) are elementary.
In fact, the lower estimate in (\ref{5.2}) has already been given in (\ref{2.6});
alternatively the upper estimate in (\ref{5.1}) follows from the Schwarz lemma.

It remains to prove the upper estimate on $e^{-u(p)}$ in (\ref{5.2}),
and we turn to that task.  Note that a dilation of $\CC$ multiplies
all quantities in (\ref{5.2}) by the same factor, so there is no loss 
of generality in assuming $\text{dist}(p,\bW)=1/2$.  Say $q_1\in\bW,\ 
|p-q_1|=1/2$.  Note that $\CC\setminus \Omega$ is connected and not bounded, 
so there exists $q_2\in \CC\setminus\Omega$ such that $|q_1-q_2|=1$.
Now, as noted in the proof of Proposition \ref{p4.1}, we have
$e^{2u}\ge e^{2w}$ on $\Omega$, where $e^{2w}\delta_{jk}$ is the 
Poincar{\'e}  metric on $\CC\setminus\{q_1,q_2\}$.  In view of the obvious
relation between the Poincar{\'e} metric on $\CC\setminus\{q_1,q_2\}$
and the Poincar{\'e} metric (call it $\Phi^2\delta_{jk}$) on $\CC\setminus
\{0,1\}$, we have
\begin{equation}
\text{dist}(p,\bW)=\frac{1}{2}\Longrightarrow e^{u(p)}\ge 
\inf\limits_{|z|=1/2}\, \Phi(z)>0,
\label{5.5}
\end{equation}
and the proof is complete.

$\text{}$\newline
{\bf Remark.} The simple argument above reveals the role of 
the simple connectivity of $\Omega$ in the estimate on the 
Poincar{\'e} metric.  In \cite{BP} there is a more sophisticated
argument yielding an estimate on the Poincar{\'e} metric of
planar domains that are not simply connected.
\newline $\text{}$ 

One classical use of Koebe's disk theorem is to provide a 
uniform bound on $|f(z)|$ for $f\in{\cal{S}}$.  We record 
a short derivation of such a bound from (\ref{5.4}).  To begin, 
the second inequality in (\ref{5.4}) implies $(1-|z|^2)|f'(z)|
\le b(1+|f(z)|)$.  Noting that $f(D_{1/2})\supset D_{1/2b}$,
we see that, for all $f\in{\cal{S}}$, 
\begin{equation}
|z|\ge \frac{1}{2}\Rightarrow 
|f(z)|\ge \frac{1}{2b}
\Rightarrow\Bigl|\frac{f'(z)}{f(z)}\Bigr|
\le \frac{B}{1-|z|^2},\quad B=b(1+2b).
\label{5.6}
\end{equation}
Since $f(D_{1/2})$ does not contain $\overline{D}_{1/2}$, we see
that there exists $z(f)$ such that $|z(f)|=1/2$ and $|f(z(f))|=1/2$.
Using (\ref{5.6}) and integrating over an arc of $\{z:|z|=1/2\}$,
we have an absolute bound
\begin{equation}
|z|=\frac{1}{2}\Longrightarrow |f(z)|\le C,\quad \forall\ f\in{\cal{S}}.
\label{5.7}
\end{equation}
This bound also holds for $|z|\le 1/2$.  Then radial integration of
(\ref{5.6}) yields an absolute bound
\begin{equation}
|f(z)|\le C_2(1-|z|)^{-B/2},
\quad \forall\ f\in{\cal{S}}.
\label{5.8}
\end{equation}

There exist sharper bounds on elements of ${\cal{S}}$, obtained by
harder work; cf.~\cite{Ahl2}, p.~84.  However, the bound derived
above suffices for the following normal family result, also due
to Koebe.  We record the essentially standard proof.

\begin{proposition} The set ${\cal{S}}$ is compact in ${\cal{H}}(D_1)$,
the space of holomorphic functions on $D_1$.
\label{p5.4}
\end{proposition}

\noindent
{\bf Proof:} Take $f_\nu\in{\cal{S}}$.  The uniform bounds $|f_\nu(z)|\le 
K(r)$ for $|z|\le r$ established above imply some subsequence converges
locally uniformly to $f\in{\cal{H}}(D_1)$.  We have $f(0)=0$ and $f'(0)=1$.
That $f$ is univalent is then a simple consequence of Hurwitz' theorem.

$\text{}$ \newline
{\bf Remark.} A direct proof of Proposition \ref{p5.4}, not using 
Proposition \ref{p5.1}, but somewhat longer and trickier, is given
in \cite{F}.

\section{The uniformization theorem}

In this section we prove the uniformization theorem for general noncompact
Riemann surfaces:

\begin{theorem} If $M$ is a noncompact, simply connected Riemann
surface, then $M$ is holomorphically equivalent to either $D_1$ or $\CC$.
\label{t6.1}
\end{theorem}

To prove this, we begin by taking $\Omega_\nu\subset\subset\Omega_{\nu+1}
\nearrow M$, such that each set $\Ombar_{\nu}$ is compact, with smooth 
boundary, and simply connected.  (This relies on some results on the topology
of surfaces, such as the a priori knowledge that $M$ is diffeomorphic to 
$D_1$).  Our argument from here parallels one in \cite{F}, except that we apply 
the method of the curvature equation to each $\Ombar_{\nu}$.

In detail, say $p\in \Omega_0\subset\subset\Omega_1\subset\subset\cdots$.
By Proposition \ref{p2.4} we have for each $\nu$ a holomorphic diffeomorphism
\begin{equation}
\psi_\nu:\Omega_\nu\longrightarrow D_1,\quad \psi_\nu(p)=0.
\label{6.1}
\end{equation}
Take $\alpha_\nu=D\psi_\nu(p)\in \mbox{Hom}(T_pM,\CC)$.  Then $\alpha_\nu=
a_\nu \alpha_0$ for uniquely defined $a_\nu\in\CC$, and if we set
\begin{equation}
\varphi_\nu:\Omega_\nu\longrightarrow D_{R_\nu},\quad
R_\nu=|a_\nu|^{-1},\quad \varphi_\nu(x)=a_\nu^{-1} \psi_\nu(x),
\label{6.2}
\end{equation}
we have
\begin{equation}
D\varphi_\nu(p)=D\varphi_0(p)=\alpha_0,\quad \forall\ \nu.
\label{6.3}
\end{equation}
It follows from the Schwarz lemma that $|a_0|>|a_1|>|a_2|>\cdots$, and hence
\begin{equation}
R_0<R_1<R_2<\cdots.
\label{6.4}
\end{equation}

To proceed, let us consider
\begin{equation}
\Phi_\nu=\varphi_\nu\circ\varphi_0^{-1}:D_1\longrightarrow D_{R_\nu}.
\label{6.5}
\end{equation}
We have each $\Phi_\nu$ holomorphic and one-to-one (i.e., univalent), and
\begin{equation}
\Phi_\nu(0)=0,\quad \Phi'_\nu(0)=1.
\label{6.6}
\end{equation}
At this point we apply Koebe's normal family theorem, established 
in the last section.  Thus we see that a subsequence of
$\varphi_\nu\bigr|_{\Omega_0}$ converges to a univalent map $\Omega_0
\rightarrow\CC$.  A similar consideration of 
\begin{equation}
\Phi_{\nu\mu}=\varphi_\nu \circ \varphi_\mu^{-1}:D_{R_\mu}\longrightarrow
D_{R_\nu},\quad \nu\ge\mu,
\label{6.7}
\end{equation}
plus a diagonal argument yields a subsequence of $(\varphi_\nu)$ 
converging to a univalent holomorphic map
\begin{equation}
\varphi:M\longrightarrow \CC.
\label{6.8}
\end{equation}

From here one could argue that $\varphi$ maps $M$ biholomorphically onto
$D_R$, where $R$ is the supremum of the sequence (\ref{6.4}); cf.~\cite{F}.
For our purposes we can bypass this argument. At this point we have $M$
holomorphically equivalent to $\Omega=\varphi(M)\subset \CC$, and the results
of \S{4} are applicable, to show that either $\Omega=\CC$ or $\Omega$ is
holomorphically equivalent to $D_1$.

\section{The curvature dichotomy}

Let $\Omega$ be a noncompact Riemann surface, with a compatible
Riemannian metric $g_0$.
Take compact, smoothly bounded $\Ombar_\nu\subset\Omega$ such that
$\Omega_\nu\nearrow\Omega$, and let $u_\nu\in C^\infty(\Omega_\nu)$ be 
solutions to the curvature equation (\ref{1.2}) such that $e^{2u_\nu} g_0$ is a 
Poincar{\'e} metric on $\Omega_\nu$.  As we have seen, $u_\nu\searrow$
as $\nu\nearrow\infty$.  We tackle the issue of convergence of $u_\nu$.

\begin{proposition} For each $\Omega$ one of the following must
happen:
\begin{itemize}
\item[1)] $u_\nu\searrow u\in C^\infty(\Omega)$, where $u$ satisfies 
(\ref{1.2}), or 
\item[2)]  $u_\nu\searrow-\infty$ on $\Omega$.
\end{itemize}
In case 2), $\Omega$ has no metric conformal to $g_0$ with Gauss curvature
$\le -1$.
\label{p7.1}
\end{proposition}
\noindent {\bf Proof:} We already know that if all $u_\nu\ge v$ with $v$ locally
bounded then case (1) holds.
Suppose $u_\nu(p_\nu)\rightarrow -\infty$ for some sequence  
$p_\nu\in\mathcal{O}\subset\subset\Omega_N$; from now on we take 
$\nu\ge N+1$.  Since $u_\nu\le u_{N+1}$ for $\nu\ge N+1$, we have a
uniform upper bound on $\Ombar_N$: $u_\nu\le A_N<\infty$ for $\nu\ge N+1$.
We thus have a bound
\begin{equation}
|e^{2u_\nu}+k|\le A_{2N}\ \mbox{ on }\ \Ombar_N,\quad \nu\ge N+1.
\label{7.1}
\end{equation}
Hence we can find $v_\nu\in C^1(\Ombar_N)$, for $\nu\ge N+1$, such that
\begin{equation}
\Delta v_\nu=e^{2u_\nu}+k\ \mbox{ on }\ \Ombar_N,\quad 
\|v_\nu\|_{L^\infty(\Omega_N)}\le A_{3N}.
\label{7.2}
\end{equation}
Hence
\begin{equation}
\Delta(u_\nu-v_\nu)=0\ \mbox{ on }\ \Omega_N,\quad
-\infty<u_\nu-v_\nu\le A_{4N},
\label{7.3}
\end{equation}
for $\nu\ge N+1$.  That $u_\nu(p_\nu)\rightarrow -\infty$ implies $u_\nu
\rightarrow -\infty$ on $\overline{\mathcal{O}}$ now follows from Harnack's
estimate.

As for the last assertion of Proposition \ref{p7.1}, note that if $e^{2w}g_0$ has 
Gauss curvature $\le -1$ on $\Omega$ then $u_\nu\ge w$ for all $\nu$.
\hfill $\Box$

Making use of Theorem \ref{t6.1}, we can restate the dichotomy in 
Proposition \ref{p7.1} as follows.

\begin{proposition} In the setting of Proposition \ref{p7.1}, case (1) holds 
if and only if $D_1$ covers $\Omega$ and case (2) holds if and only if
$\CC$ covers $\Omega$.
\label{p7.2}
\end{proposition}

$\text{}$ \newline
{\bf Proof:} Suppose case (1) holds; we show $\Omega$ cannot be covered by 
$\CC$.  Indeed, a holomorphic covering $f:\CC\rightarrow \Omega$ would
pull back the metric $e^{2u}\, g_0$ to a metric on $\CC$ of curvature
$-1$.  However, as one sees by looking at the metrics (\ref{2.7}) on the disks
$D_\beta$ and letting $\beta\nearrow\infty$, case (2) holds for $\CC$,
so this is not possible.

On the other hand, suppose case (2) holds; we claim $\Omega$ cannot be
covered by $D_1$.  Indeed, a holomorphic covering $f:D_1\rightarrow\Omega$
puts a Poincar{\'e} metric on $\Omega$, which we have seen cannot hold
in case (2).

Since case (1) yields a holomorphic covering $f:D_1\rightarrow\Omega$, 
which puts a Poincar{\'e} metric $e^{2v}\, g_0$ on $\Omega$, and since 
$u_\nu\ge v$ on $\Omega_\nu$, we have $u\ge v$, and hence $u=v$, so:

\begin{corollary}  In case (1), the limit $u$ yields a Poincar{\'e} 
metric $e^{2u}\, g_0$ on $\Omega$.  In particular, whenever $\Omega$
has a metric of curvature $\le -1$, it has a conformally equivalent 
Poincar{\'e} metric.
\label{c7.3}
\end{corollary}

When $\Omega$ is simply connected, the dichotomy in Proposition \ref{p7.1} is 
precisely the dichotomy between hyperbolic and parabolic Riemann surfaces,
defined as follows.

$\text{}$ \newline
{\bf Definition.} A noncompact (connected) Riemann surface $\Omega$
is hyperbolic if and only if there exists a nonconstant, nonpositive, 
subharmonic function on $\Omega$, and parabolic otherwise.

$\text{}$\newline
The proof of the uniformization theorem found in most sources (for example,
\cite{Ahl2}, \cite{FK}, \cite{Sp}) proceeds by separately treating these
two cases.  It is shown that a simply connected hyperbolic Riemann surface
is equivalent to $D_1$ and a simply connected parabolic Riemann surface 
is equivalent to $\CC$.  In these treatments, the first step in the 
hyperbolic case is the construction of a negative Green function $u$,
harmonic on $\Omega\setminus\{p\}$ and behaving like $\log |x|$ in local
normal coordinates centered at $p$.  Such a function has a harmonic
conjugate $v$, well defined mod $2\pi\ZZ$, on $\Omega\setminus\{p\}$,
and the function
\begin{equation}
f(x)=e^{u(x)+iv(x)}
\end{equation}
yields a holomorphic map $f:\Omega\rightarrow D_1$.  It is then shown that
this map is a holomorphic diffeomorphism.  This argument is highly 
nontrivial, much more subtle than the demonstration of the analogue 
in the context of (\ref{2.17}). However, \cite{H} produced an ingenious 
demonstration, used in most modern treatments.  A short argument
to prove the parabolic case is given in \cite{Oss}, and other proofs
can be found in the references cited above. 
By contrast, the proof given in \S{6} of this paper is in some respects
closer to Koebe's original demonstration; compare the treatment in
\cite{Tsu}, pp.~421--422.

\section{Compact Riemann surfaces}

Here we discuss the uniformization theorem for compact Riemann surfaces,
given by the following classical result.

\begin{proposition} Let $M$ be a compact, connected Riemann surface, 
with Euler characteristic $\chi(M)$.

\begin{itemize}
\item[1)] If $\chi(M)=2$, then $M$ is holomorphically equivalent to the 
Riemann sphere.

\item[2)] If $\chi(M)=0$, then $M$ is holomorphically equivalent to a flat
torus $\TT_\Lambda=\CC/\Lambda$.

\item[3)] If $\chi(M)<0$, then $M$ is holomorphically covered by $D_1$.
\end{itemize}
\label{p8.1}
\end{proposition}

In cases 2) and 3), the universal cover $\widetilde{M}$ of $M$ is noncompact
and Theorem \ref{t6.1} applies.  If we allow ourselves to use the topological
classification of surfaces, we know that in case 3) the group of covering
transformations of $\widetilde{M}\rightarrow M$ is noncommutative, so
there is a noncommutative discrete group of holomorphic automorphisms of 
$\widetilde{M}$, acting with no fixed points.  A holomorphic automorphism
of $\CC$ without fixed points must be a translation, and all translations
commute, so in this case $\widetilde{M}$ must be conformally equivalent to 
$D_1$.  On the other hand, in case 2) the group of covering transformations 
of $\widetilde{M}\rightarrow M$ is isomorphic to $\ZZ^2$.  One can check that 
there is no group of fixed-point-free holomorphic automorphisms of $D_1$
isomorphic to $\ZZ^2$.

Case 1) is often proved via the Riemann-Roch theorem, which implies that,
given $p\in M$, the space of meromorphic functions with at most one simple 
pole, at $p$, has dimension $2$, provided $\chi(M)=2$.  This space includes 
the constant functions, but it must also contain a meromorphic function $f$
on $M$ with one simple pole (at $p$).  Then $f$ defines a holomorphic map 
$f:M\rightarrow \widehat{\CC}$.  We see that $f$ has degree one, and it 
follows that $f$ is a holomorphic diffeomorphism.

In fact, all cases of Proposition \ref{p8.1} have PDE proofs where one starts with 
some Riemannian metric on $M$ compatible with the given conformal structure.
We mention PDE treatments of these three cases.

To treat case 1), we note that given a distribution $\delta'$ of order $1$ 
supported at $p\in M$ such that $\langle 1,\delta'\rangle=0$ (a derivative
of a delta function) we can solve $\Delta u=\delta'$.  This can be done on 
any compact, connected $M$, but in case 1) we can say that
$M\setminus\{p\}$ is simply connected.
Hence $u$ is the real part of a meromorphic function $f$ on $M$ with one 
simple pole (at $p$), and as noted in our previous discussion of case 1),
this yields a holomorphic diffeomorphism $f:M\rightarrow\widehat{\CC}$.

To treat case 2), we can solve the curvature equation, which in this case
is $\Delta u=K_0(x)$.  Note that by the Gauss-Bonnet theorem $\int_M
K_0\,dA=0$, so we can solve for $u$.  Hence $M$ has a conformally equivalent
flat metric, from which the conclusion follows.

Case 3) can be established by solving the curvature equation (\ref{1.2}).  
This was accomplished by \cite{B}; expositions can also be found in 
\cite{Au} and in Chapter 14 of \cite{T}.  
Here we merely mention that the solution to 
the curvature equation (\ref{1.2}) minimizes
$$
F(u)=\int\limits_M \Bigl(\frac{1}{2}|du|^2 +K_0(x)u\Bigr) dA
$$
on the set
$$
S=\Bigl\{u\in H^1(M):\int\limits_M e^{2u}\, dA=-2\pi\chi(M)\Bigr\}.
$$

There are other very interesting analytic avenues to uniformization in the
compact case; these include the use of the determinant of the Laplacian
\cite{OPS} and a fourth-order heat flow \cite{Ch}.

\end{document}